\documentclass{article}
\usepackage{amsmath}
\usepackage[dvipsnames]{xcolor}
\usepackage{hyperref}
\usepackage{graphicx} 
\usepackage{booktabs}
\usepackage{multirow}
\usepackage{tabularx}
\usepackage{mathtools,bm}
\usepackage{algpseudocode}
\usepackage{amsthm,amsfonts,amssymb}
\usepackage{algorithm}

\newtheorem{theorem}{Theorem}
\newtheorem{proposition}[theorem]{Proposition}%
\newtheorem{corollary}[theorem]{Corollary}
\newtheorem{remark}{Remark}%

\numberwithin{equation}{section}

\title{A fast surrogate cross validation algorithm\\ for meshfree RBF collocation approaches}
\author{Francesco Marchetti\\
Department of Mathematics \lq\lq Tullio Levi-Civita\rq\rq, Università di Padova \thanks{francesco.marchetti@unipd.it}}
\date{}

\begin{document}

\maketitle


\section*{Abstract}
{Cross validation is an important tool in the RBF collocation setting, especially for the crucial tuning of the shape parameter related to the radial basis function. In this paper, we define a new efficient surrogate cross validation algorithm, which computes an accurate approximation of the true validation error with much less computational effort with respect to a standard implementation. The proposed scheme is first analyzed and described in details, and then tested in various numerical experiments that confirm its efficiency and effectiveness.}


\section{Introduction}

Radial Basis Function (RBF) collocation methods are widely-investigated approaches for solving various Partial Differential Equations (PDEs). They are kernel-based \textit{meshfree} methods, flexible and potentially highly accurate, which deal with PDEs in their \textit{strong form} by \textit{interpolating} at a set of collocation points. For this reason, they can be intended and studied in a generalized interpolation framework \cite{Fasshauer07,Fasshauer15}. 

For both collocation and interpolation tasks, a crucial aspect for the construction of an effective approximant is an appropriate choice of the RBF, whose form and behavior very often depend on a positive shape parameter. Although a \textit{variable} shape parameter is employed in some cases \cite{Chen22a,Chiu20,Kansa92}, in literature greater attention has been put on the selection of a \textit{global} value, and many strategies have been proposed for its fine tuning  \cite{Chen22,Fornberg07,Gherlone12,Krowiak19,Scheuerer11,Yang18}, including some theoretical observations and \textit{rules of thumb} \cite{Golbabai15,Uddin14}. In particular, many of these strategies take advantage of Cross Validation (CV) approaches \cite{Golub79}, which are empirical schemes that are widely employed in numerous fields in order to estimate the generalization capability of a model. Historically, in the RBF interpolation setting, the most considered CV approach is the Leave-One-Out CV (LOOCV) method, mainly because of its efficient computation provided by an algorithm originally proposed by S. Rippa in \cite{Rippa99}, and modified by some authors later on \cite{Mongillo11} . In \cite{Marchetti21}, the Extended Rippa's Algorithm (ERA) generalized Rippa's scheme to the wider $k$-fold CV framework, still including in its formulation the original LOOCV case.

Because of its efficiency, Rippa's approach has been extensively employed beyond the interpolation setting, with the idea of obtaining an \textit{inexact} but effective CV tool. In \cite{Cavoretto22,Cavoretto20a,Cavoretto20b}, the interpolation matrix used in Rippa's LOOCV scheme is directly substituted by the collocation matrix, with applications to elliptic and time-dependent PDEs, also on irregular domains. Many other works share a similar spirit, \cite{Chen20,Roque10,Trahan03} are some examples. In \cite{Fasshauer07a} the authors adapted Rippa's approach to obtain a cost indicator in the context of RBF-PseudoSpectral (RBF-PS) methods and approximate moving least squares, and a modified scheme was designed in the partition of unity interpolation framework in \cite{Cavoretto18}. However, a concrete and \textit{justified} generalization of the ERA to the collocation setting is still missing.

In this paper, our purpose is to fill this gap and to provide a generalization of the ERA scheme for PDEs, by analyzing to what extent the original algorithm can be adapted from RBF interpolation to collocation. The result is a surrogate CV scheme that retains the efficiency of the ERA, but computes an approximated validation error. Nevertheless, differently from the adaptations already suggested in literature, the discrepancy between the validation error and its approximation can be clearly formalized. In fact, maybe quite surprisingly, we show that this \textit{approximation gap} involves the reconstruction provided by the RBF-PS method. Various numerical experiments show that the proposed technique well approximates exact CV, being more accurate than the empirical adaptation of Rippa's approach proposed for the collocation framework in literature. Moreover, in some situations our scheme can also be employed when collocation points do not coincide with centers.

The paper is organized as follows. In Section \ref{sec:int_col} we recall the main characteristics of the generalized RBF interpolation setting, and how classical collocation approaches can be framed into it. In Section \ref{sec:gerameio} we present and analyze in details the proposed \textit{Surrogate CV} scheme, both from a theoretical and a computational perspective. Various numerical tests are discussed in Section \ref{sec:numerics}. Finally, in Section \ref{sec:conclusions} we draw some conclusions.

\section{Generalized interpolation and collocation}\label{sec:int_col}
Let $\Omega\subset \mathbb{R}^d$ and let $\kappa:\Omega\times\Omega\longrightarrow\mathbb{R}$ be a strictly positive definite kernel that is radial, i.e., there exists a univariate Radial Basis Function (RBF) $\varphi:[0,\infty)\longrightarrow \mathbb{R}$ so that 
$$\kappa(\bar{x},\bar{y})=\varphi(\lVert \bar{x}-\bar{y} \lVert ),\qquad \bar{x},\bar{y}\in\Omega.$$
The RBF $\varphi$ usually depends on a positive shape parameter $\varepsilon>0$. However, to keep a simpler notation, we avoid to make this dependency explicit until Section \ref{sec:numerics}, where $\varepsilon$ will be of interest in our numerical experiments. We can introduce the framework of kernel-based collocation as a generalized interpolation setting. Let $\mathcal{G}=\{\gamma_1,\dots,\gamma_m\}$ and $\mathcal{L}=\{\lambda_1,\dots,\lambda_n\}$ be sets of distinct linear functionals acting on real-valued functions defined on $\Omega$, $m,n\in\mathbb{N}$. Given $h:\Omega\longrightarrow \mathbb{R}$, our purpose is to find $s:\Omega\longrightarrow \mathbb{R}$ so that
\begin{equation}\label{eq:collocation}
	\gamma_i(h)=\gamma_i(s)\qquad i=1,\dots,m.
\end{equation}
If we restrict to the kernel-based \textit{ansatz}
\begin{equation}\label{eq:ansatz}
	s(\bar{x})=\sum_{j=1}^n c_j\lambda_j(\kappa(\bar{x},\bar{y})),
\end{equation}
where $\lambda_j$ acts on $\kappa$ seen as a function of its second input and $c_j\in\mathbb{R}$, then \eqref{eq:collocation} becomes
\begin{equation*}
	\gamma_i(h)=\sum_{j=1}^n c_j(\gamma_i\circ\lambda_j)(\kappa(\bar{x},\bar{y})),\qquad i=1,\dots,m,
\end{equation*}
with $\gamma_i$ acting on $\kappa$ seen as a function of its first input. Therefore, the vector of coefficients $\bar{c}=(c_1,\dots,c_n)^\intercal$ needs to satisfy
\begin{equation}\label{eq:lin_col}
	\mathsf{G}\bar{c}=\bar{g},
\end{equation}
being $\mathsf{G}_{i,j}=(\gamma_i\circ\lambda_j)(\kappa(\bar{x},\bar{y}))$ the \textit{collocation} matrix and $\bar{g}=(\gamma_1(h),\dots,\gamma_m(h))^\intercal$. 

In this context, the considered functionals are usually related to sets of points in $\Omega$. We then associate the set of collocation points $X=\{\bar{x}_1,\dots,\bar{x}_m\}\subset\Omega$ and \textit{centers} $Y=\{\bar{y}_1,\dots,\bar{y}_n\}\subset\Omega$ to $\mathcal{G}$ and $\mathcal{L}$, respectively. More precisely, $\gamma_i=\gamma_{\bar{x}_i}$ and $\lambda_j=\lambda_{\bar{y}_j}$ for each $i=1,\dots,m$ and $j=1,\dots,n$. We point out that $\mathcal{G}$ and $\mathcal{L}$ may be related to the same set of points, even if they do not coincide. 

Another relevant tool in our analysis is then the \textit{evaluation} matrix $\mathsf{L}_{i,j}=(\delta_{\bar{x}_i}\circ\lambda_{\bar{y}_j})(\kappa(\bar{x},\bar{y}))$, where $\delta_{\bar{x}_i}(f)=f(\bar{x}_i)$, $i=1,\dots,m$, are point evaluation functionals at the collocation points in $X$. By defining $\bar{s}=(s(\bar{x}_1),\dots,s(\bar{x}_m))^\intercal$, note that we can write
\begin{equation}\label{eq:lin_eva}
	\bar{s}=\mathsf{L}\bar{c}.
\end{equation}
In general, there is no guarantee that \eqref{eq:lin_col} admits one or more solutions, and therefore the representation of $\bar{s}$ provided in \eqref{eq:lin_eva} might not exist, or it might be non unique.

Let us outline two well-known collocation approaches that fall into the presented framework. We will consider them in our numerical experiments in Section \ref{sec:numerics}.

\subsection{Kansa's method: \lq\lq$\gamma\neq \lambda,\; \lambda=\delta$\rq\rq}\label{sec:kansa}

In the approach originally proposed by Ed Kansa in \cite{Kansa90}, the functionals in $\mathcal{L}$ are chosen as $\lambda_{\bar{y}_j}=\delta_{\bar{y}_j}$, $j=1,\dots,n$. The model \eqref{eq:ansatz} then takes the form
\begin{equation*}
	s(\bar{x})=\sum_{j=1}^n c_j\kappa(\bar{x},\bar{x}_j),
\end{equation*}
which is the same of function interpolation and regression settings. Therefore, this approach takes advantage of a simple ansatz. On the other hand, since the functionals in $\mathcal{G}$ are employed to apply the differential operators that characterize the PDE, the resulting collocation matrix $\mathsf{G}$ is typically non-symmetric (see e.g. \cite[Section 38.1]{Fasshauer07}). Even if it has been shown that is unlikely for $\mathsf{G}$ to not have full rank \cite{Hon01}, this represented a relevant obstacle for the theoretical analysis of this approach \cite{Schaback07}. Asides from these issues, Kansa's approach has been tailored to various differential problems and largely employed by the relevant scientific community \cite{Chiappa20,Dehghan15,Dehghan14,Karageorghis21,Liu15,Wang18}.

\subsection{Hermite method: \lq\lq$\gamma= \lambda\neq \delta$\rq\rq} 

In the Hermite approach, both the functionals in $\mathcal{G}$ and in $\mathcal{L}$ are used to represent the operators involved in the PDE (see, e.g., \cite[Section 38.2]{Fasshauer07}). A great advantage of this method is that the collocation matrix is always symmetric and invertible in the case where collocation points and centers coincide, and therefore \eqref{eq:lin_col} is well posed (refer to e.g. \cite{Wendland05}). On the other hand, applying differential operators twice on the first and second input of the kernel may lead to a solution that is double regular than needed, and this can represent a drawback. Similarly to Kansa's method, the Hermite approach caught the interest of many authors in the last years \cite{Chu14,LaRocca05,LaRocca06,Kazem12,Ma21}.

\section{Efficient surrogate cross validation}\label{sec:gerameio}

\subsection{Case $X=Y$}\label{sec:xy}

It is convenient to start from the case where collocation points and centers coincide, i.e., $X=Y$ and $m=n$. In this setting both $\mathsf{G}$ and $\mathsf{L}$ are square matrices, and $\mathsf{G}$ is symmetric in the case of the Hermite method.

Before presenting in details the proposed algorithm, we recall the main characteristics of classical Cross Validation (CV). Let $k\in\mathbb{N}$, $k\le m$. In the $k$-fold CV setting, the set $X$ is split into $k$ disjoint non-empty subsets $X_1,\dots,X_k$. For a simpler presentation, we assume $m$ to be multiple of $k$, so that we can consider subsets of equal cardinality $v=m/k\in\mathbb{N}$, but this is not required by the scheme. Then, for every $\ell=1,\dots,k$ we proceed as follows.
\begin{enumerate}
	\item 
	We construct the approximant $s$ by employing only the nodes in $\bigcup_{\substack{i=1\\ i\neq\ell}}^k X_\ell$.
	\item 
	We assess the validation performance of $s$ on the remaining set $X_\ell$ by computing, e.g., the $L_2$-error with respect to the \textit{ground truth} $h$.
\end{enumerate}
At the end of this process, we computed and evaluated $k$ independent models exploiting the information in $X$ only. The \textit{global} performance on $X$ obtained by putting together the $k$ validation errors is then employed as an indicator of the overall accuracy of the model when constructed on the full set of collocation points $X$. In this sense, $k$-fold CV is a \textit{fair} and concrete procedure that does not assume knowledge at unknown data sites, differently with respect to, e.g., \textit{trial-and-error} approaches.

The computational complexity required by the whole process is $\mathcal{O}(k(m-v)^3)\approx\mathcal{O}(m^3k)$: we need to solve $k$ different $(m-v)\times (m-v)$ linear systems (we assume the usage of non-specific techniques such as LU decomposition to solve the linear systems). In RBF applications it is very common to set $v\ll m$, and in particular when $v=1$ we get the so-called Leave-One-Out CV (LOOCV). Consequently, the computational cost required by a straightforward application of LOOCV is approximately $\mathcal{O}(m^4)$.

In the framework of kernel-based interpolation, the Extended Rippa's Algorithm (ERA) computes exact $k$-fold CV at the cost of $\mathcal{O}(m^3)+\mathcal{O}\big(\frac{m^3}{k^2}\big)$ operations, which is indeed really advantageous especially when $v$ is small.
Inspired by the ERA, our purpose is to construct an algorithm for validating RBF collocation schemes and enhancing the efficiency of a standard $k$-fold CV approach. 

Let us consider a single validation step of the CV procedure (a certain $\ell\in\{1,\dots,m\}$), and let us fix some notations. We define as $\bar{p}=(p_1,\dots,p_v)^{\intercal},\;p_i\in\{1,\dots,n\}$ the vector of distinct validation indices that identifies the elements of the validation set $X_\ell=\{\bar{x}_{p_1},\dots,\bar{x}_{p_v}\}$. Furthermore, given a $m$-dimensional vector $\bar{z}$ and a $m\times m$ matrix $\mathsf{A}$, we denote as:
\begin{itemize}
	\item 
	$\bar{z}_{\bar{p}}$ the $v$-dimensional vector obtained by restricting to the elements whose index is in $\bar{p}$, and $\bar{z}^{\bar{p}}$ the $(m-v)$-dimensional vector obtained by restricting to the elements whose index is not in $\bar{p}$.
	\item 
	$\mathsf{A}_{\bar{p},\bar{p}}$ the $v\times v$ matrix obtained by restricting to the rows and columns whose index is in $\bar{p}$, and $\mathsf{A}^{\bar{p},\bar{p}}$ the $(m-v)\times (m-v)$ matrix obtained by restricting to the rows and columns whose index is not in $\bar{p}$.
\end{itemize}
The introduced notation is helpful for presenting the following result. 
\begin{theorem}\label{thm:core}
	Assume that the $m\times m$ matrices $\mathsf{G}$ and $\mathsf{L}$ in \eqref{eq:lin_col} and \eqref{eq:lin_eva} are invertible. Let $s$ be the approximant \eqref{eq:ansatz} constructed at the full set of collocation points/centers $X$ by solving the linear system \eqref{eq:lin_col}, and let $s^{(\bar{p})}$ be the approximant built by excluding the functionals related to the nodes in $X_\ell$, i.e.,
	\begin{equation*}
		s^{(\bar{p})}(\bar{x})=\sum_{\bar{x}_j\notin X_\ell} c_j^{(\bar{p})}\lambda_{\bar{x}_{j}}(\kappa(\bar{x},\bar{y})),
	\end{equation*}
	where the vector of coefficients $\bar{c}^{(\bar{p})}=(c_1^{(\bar{p})},\dots,c_v^{(\bar{p})})^\intercal$ solves the linear system $\mathsf{G}^{\bar{p},\bar{p}}\bar{c}^{(\bar{p})}=\bar{g}^{\bar{p}}$. Moreover, let $\bar{h}=(\delta_1(h),\dots,\delta_m(h))^\intercal$ be the vector of evaluations related to the underlying function $h$. Then, the vector of signed validation errors $\bar{e}_{\bar{p}}=\bar{s}^{(\bar{p})}(X_\ell)-\bar{g}_{\bar{p}}$, with $\bar{s}^{(\bar{p})}(X_\ell)=(\bar{s}^{(\bar{p})}(\bar{x}_{p_1}),\dots,\bar{s}^{(\bar{p})}(\bar{x}_{p_v}))^\intercal$, can be approximated as
	\begin{equation}
		\bar{e}_{\bar{p}}\approx\bar{\epsilon}_{\bar{p}}=((\mathsf{L}^{-1})_{:,\bar{p}})^+(\mathsf{G}^{-1})_{:,\bar{p}}((\mathsf{G}^{-1})_{\bar{p},\bar{p}})^{-1}(\mathsf{G}^{-1}\bar{g})_{\bar{p}}-\bar{f}_{\bar{p}}+\bar{h}_{\bar{p}},
	\end{equation}
	where $\mathsf{A}^+$ denotes the Moore-Penrose inverse (or pseudoinverse) of the matrix and $\bar{f}=\mathsf{L}\mathsf{G}^{-1}\bar{g}$ is the RBF-PS solution at $X$. Precisely, $\bar{e}_{\bar{p}}=\bar{\epsilon}_{\bar{p}}$ if the residual
	\begin{equation*}
		\lVert \big((\mathsf{L}^{-1})_{:,\bar{p}}((\mathsf{L}^{-1})_{:,\bar{p}})^+-\mathsf{I}\big)(\mathsf{G}^{-1})_{:,\bar{p}}((\mathsf{G}^{-1})_{\bar{p},\bar{p}})^{-1}(\mathsf{G}^{-1}\bar{g})_{\bar{p}}\lVert_2
	\end{equation*}
	is equal to zero, being $\mathsf{I}$ the $m\times m$ identity matrix.
\end{theorem}

\begin{proof}
	Let $\bar{f}=\mathsf{L}\bar{c}\in\mathbb{R}^m$, which implies $\bar{f}= \mathsf{L}\mathsf{G}^{-1}\bar{g}$. Moreover, let $\bar{b}=(b_1,\dots,b_m)^{\intercal}\in\mathbb{R}^m$ be such that:
	\begin{enumerate}
		\item[(P1)] $\bar{b}_{\bar{p}} \equiv \bar{0}$.
		\item[(P2)] $\mathsf{L}\bar{b} =  \bar{f} - \sum_{j=1}^v{\alpha_j \mathsf{I}_{:,p_j}}$, where $\mathsf{I}_{:,p_j}$ denotes the $p_j$-th column of the $m\times m$ identity matrix $\mathsf{I}$, and $\bar{\alpha}=(\alpha_1,\dots,\alpha_v)^{\intercal}\in\mathbb{R}^v$.
		\item[(P3)] $\mathsf{G}\bar{b} =  \bar{g} - \sum_{j=1}^v{\beta_j \mathsf{I}_{:,p_j}}$, where $\bar{\beta}=(\beta_1,\dots,\beta_v)^{\intercal}\in\mathbb{R}^v$.
	\end{enumerate}
	 Note that for any $m\times m$ matrix $\mathsf{A}$ and vector $\bar{z}$ such that $\mathsf{A}\bar{b}=\bar{z}$, (P1) implies $\mathsf{A}^{\bar{p},\bar{p}}\bar{b}^{\bar{p}}=\bar{z}^{\bar{p}}$. By taking into account also (P2) and recalling \eqref{eq:lin_eva}, we obtain
	 \begin{equation*}
	 	s^{(\bar{p})}(\bar{x}_{p_i})=\sum_{\bar{x}_j\notin X_\ell} c_j^{(\bar{p})}(\delta_{\bar{x}_{p_i}}\circ\lambda_{\bar{x}_{j}})(\kappa(\bar{x},\bar{y}))=\sum_{j=1}^m b_j	(\delta_{\bar{x}_{p_i}}\circ\lambda_{\bar{x}_{j}})(\kappa(\bar{x},\bar{y}))=(\mathsf{L}\bar{b})_{p_i}=\bar{f}_{p_i}-\alpha_{p_i},	
	 \end{equation*}
	 which implies
	 \begin{equation}\label{eq:alfa_bareta}
	 	\bar{\alpha}=\bar{f}_{\bar{p}}-\bar{s}^{(\bar{p})}(X_\ell)=(\mathsf{L}\mathsf{G}^{-1}\bar{g})_{\bar{p}}-\bar{s}^{(\bar{p})}(X_\ell).
	 \end{equation}
	 Now, from (P3) we get
	 \begin{equation*}
	 	\bar{b} = \mathsf{G}^{-1}\bigg(\bar{g} - \sum_{j=1}^v{\beta_j \mathsf{I}_{:,p_j}}\bigg)= \bar{c}- \sum_{j=1}^v{\beta_j(\mathsf{G}^{-1})_{:,p_j}},
	 \end{equation*}
	 and then $\bar{0}\equiv\bar{b}_{\bar{p}}= \bar{c}_{\bar{p}} - (\mathsf{G}^{-1})_{\bar{p},\bar{p}}\bar{\beta}$, from which we can calculate $\bar{\beta}=((\mathsf{G}^{-1})_{\bar{p},\bar{p}})^{-1}\bar{c}_{\bar{p}}$.
	 
	 The next step consists in analyzing $\bar{\alpha}$ in terms of $\bar{\beta}$. Therefore, we combine (P2) and (P3) and get
	 \begin{equation*}
	 	\mathsf{L}^{-1}\bigg(\bar{f} - \sum_{j=1}^v{\alpha_j \mathsf{I}_{:,p_j}}\bigg)=\mathsf{G}^{-1}\bigg(\bar{g} - \sum_{j=1}^v{\beta_j \mathsf{I}_{:,p_j}}\bigg),
	 \end{equation*}
 	from which it follows
 	\begin{equation*}
 		\sum_{j=1}^v{\alpha_j (\mathsf{L}^{-1})_{:,p_j}}=\sum_{j=1}^v{\beta_j(\mathsf{G}^{-1})_{:,p_j}}.
 	\end{equation*}
 Therefore, the sought vector $\bar{\alpha}$ has to satisfy $(\mathsf{L}^{-1})_{:,\bar{p}}\bar{\alpha}=(\mathsf{G}^{-1})_{:,\bar{p}}((\mathsf{G}^{-1})_{\bar{p},\bar{p}})^{-1}(\mathsf{G}^{-1}\bar{g})_{\bar{p}}$. This linear system is likely to be overdetermined, and to not admit an exact solution. Thus, we take the least squares solution 
 \begin{equation}\label{eq:lss}
 	\bar{\alpha}^{\star}=\min_{\bar{\alpha}}{\lVert (\mathsf{L}^{-1})_{:,\bar{p}}\bar{\alpha}-(\mathsf{G}^{-1})_{:,\bar{p}}((\mathsf{G}^{-1})_{\bar{p},\bar{p}})^{-1}(\mathsf{G}^{-1}\bar{g})_{\bar{p}}\lVert_2}.
 \end{equation}
 By employing the pseudoinverse, we can write
 \begin{equation*}
 	\bar{\alpha}^{\star}= ((\mathsf{L}^{-1})_{:,\bar{p}})^+(\mathsf{G}^{-1})_{:,\bar{p}}((\mathsf{G}^{-1})_{\bar{p},\bar{p}})^{-1}(\mathsf{G}^{-1}\bar{g})_{\bar{p}}.
 \end{equation*}
Finally, we observe that standard CV yields to
\begin{equation*}
	\bar{e}_{\bar{p}}=\bar{h}_{\bar{p}}-\bar{s}^{(\bar{p})}(X_\ell),
\end{equation*}
therefore
\begin{equation*}
	\bar{e}_{\bar{p}}-\bar{\alpha}^{\star}=\bar{h}_{\bar{p}}-\bar{f}_{\bar{p}},
\end{equation*}
and
\begin{equation*}
	\bar{\epsilon}_{\bar{p}}=\bar{\alpha}^{\star}-\bar{f}_{\bar{p}}+\bar{h}_{\bar{p}}.
\end{equation*}
\end{proof}
\begin{corollary}
	In the case of LOOCV where $\bar{p}=p\in\{1,\dots,m\}$, we get
	\begin{equation*}
		\bar{\epsilon}_{p}=((\mathsf{L}^{-1})_{:,{p}})^+(\mathsf{G}^{-1})_{:,{p}}\frac{(\mathsf{G}^{-1}\bar{g})_{{p}}}{(\mathsf{G}^{-1})_{{p},{p}}}.
	\end{equation*}
\end{corollary}
\begin{proof}
	The result is achievable by some simple manipulations, and by observing that $(\mathsf{G}^{-1})_{{p},{p}}\in\mathbb{R}$.
\end{proof}
Let us comment on the presented theorem. In general, the assumption regarding the invertibility of the collocation and evaluation matrices is not theoretically justified. Nevertheless, as also outlined in Subsection \ref{sec:kansa}, it is a common practice to assume such hypothesis in the context of RBF collocation. The Surrogate CV requires the calculation of the inverses of $\mathsf{G}$ and $\mathsf{L}$, but in each validation step $\ell=1,\dots,k$ the required computational effort is reduced with respect to the classical approach. We better analyze this important aspect in Subsection \ref{sec:comput}.

As far as the discrepancy between the error vectors computed by exact and Surrogate CV is concerned, the inexactness lies in the fact that  $\bar{\epsilon}_{\bar{p}}$ is likely to only approximate the validation error between the \textit{true} underlying function $h$ and the approximation computed by means of a RBF-PseudoSpectral (RBF-PS) approach that takes advantage of the full set $X$ (cf. e.g. \cite[Chapter 42]{Fasshauer07}). Nevertheless, we point out that the RBF-PS approximation can be calculated explicitly in our setting without adding a relevant computational cost (see Algorithm \ref{alg:eccheo}). 

\begin{remark}
	If we choose $\gamma_{\bar{x}_i}=\lambda_{\bar{x}_i}=\delta_{\bar{x}_i}$, $i=1,\dots,m$, then $\mathsf{G}=\mathsf{L}$. As a consequence, in the proof of Theorem \ref{thm:core} we get $\bar{f}=\bar{g}$ and we recover the exact ERA from the RBF interpolation framework. In this sense, the proposed Surrogate CV algorithm can be considered as a generalization of the ERA for generalized interpolation.
\end{remark}

\subsection{Case $X\neq Y$}

In the following, we investigate whether the proposed Surrogate CV scheme can be extended to the more general case where collocation points differ from centers. This setting is more difficult to treat from a formal viewpoint than the $X=Y$ case, both in interpolation and in collocation, but considering centers that do not exactly correspond to collocation points may lead to an increased numerical stability in certain situations as observed, e.g., in \cite{Campagna20,Katsiamis20}. In order to apply our scheme, the validation partition $X_1,\dots,X_v$ has to be relatable to both the rows and the columns of $\mathsf{G}$ and $\mathsf{L}$. This is clear since the first part of the proof of Theorem \ref{thm:core}, until \eqref{eq:alfa_bareta}, where the validation set $X_\ell$ and the corresponding vector of indices $\bar{p}$ have influence on both $\delta_i$ and $\lambda_j$. Therefore, we can say something in the case where the set of centers consists of the collocation points plus some other additional centers.

\begin{corollary}\label{cor:sticentri}
	Let $Y=X\cup Z$, where $Z=\{\bar{z}_1,\dots,\bar{z}_\nu\}\subset\mathbb{R}^d$, $\nu\in\mathbb{N}$, is an additional set of centers that are not considered as collocation points. Then, the Surrogate CV that relies on the collocation points set $X$ leads to
	\begin{equation*}
		\bar{\epsilon}_{\bar{p}}= ((\mathsf{L}^{+})_{:,\bar{p}})^+(\mathsf{G}^{+})_{:,\bar{p}}((\mathsf{G}^{+})_{\bar{p},\bar{p}})^{-1}(\mathsf{G}^{+}\bar{g})_{\bar{p}}-\bar{f}_{\bar{p}}+\bar{h}_{\bar{p}}.
	\end{equation*}
\end{corollary}
\begin{proof}
	The matrices $\mathsf{G}$ and $\mathsf{L}$ are now $(m\times n)$-dimensional with $n=m+\nu$. We can construct them in such a way that the last columns correspond to the functionals related to $Z$, i.e.,
	\begin{equation*}
		\begin{cases}
			\mathsf{G}_{i,j}=(\gamma_{\bar{x}_i}\circ\lambda_{\bar{x}_j})(\kappa(\bar{x},\bar{y})), & \quad $j=1,\dots,m$,\\
			\mathsf{G}_{i,j}=(\gamma_{\bar{x}_i}\circ\lambda_{\bar{z}_{j-m}})(\kappa(\bar{x},\bar{y})), & \quad $j=m+1,\dots,n$,\\
		\end{cases}
	\end{equation*}
	and doing the same $\mathsf{L}$. Consequently, we can proceed as in the proof of Theorem \ref{thm:core}, keeping in mind that:
	\begin{itemize}
		\item 
		Classical matrix inversions are replaced by pseudoinverses, with the effect that we choose minimal norm solutions among the infinite admissible ones (the related linear systems are now undetermined).
		\item 
		The validation procedure does not involve the additional centers in $Z$, that are therefore always employed as centers in every model and never excluded during the process.
	\end{itemize}
\end{proof}

\subsection{Formulating the algorithm}\label{sec:comput}

We can exploit the findings of the previous subsections to finalize the Surrogate CV method, which is detailed in Algorithm \ref{alg:eccheo}. In order to present a unique version of the scheme, we consider the more general assumptions of Corollary \ref{cor:sticentri}, observing that in the case $X=Y$ the pseudoinverses turn to classical matrix inverses when appropriate. 
\begin{algorithm}[h!]
	\caption{Surrogate $k$-fold CV}\label{alg:eccheo}
	\begin{algorithmic}[1]
		\Require $\mathsf{G}$, $\mathsf{L}$, $\bar{g}$, $\bar{h}$, $k$
		\Ensure $\bar{\epsilon}$
		\State Set $k$ vectors $\bar{p}_1,\dots,\bar{p}_k$ of distinct indices in $\{1,\dots,m\}$
		\State Compute $\mathsf{G}^+$, $\mathsf{L}^+$
		\State Compute $\bar{c}=\mathsf{G}^+\bar{g}$
		\State Compute $\bar{f}=\mathsf{L}\bar{c}$
		\For{$\ell=1,\dots,k$} 	
		\State Compute $\bar{\alpha}_{\bar{p}_\ell}= ((\mathsf{L}^{+})_{:,\bar{p}_\ell})^+(\mathsf{G}^{+})_{:,\bar{p}_\ell}((\mathsf{G}^{+})_{\bar{p}_\ell,\bar{p}_\ell})^{-1}\bar{c}_{\bar{p}_\ell}$
		\EndFor
		\State Complete $\bar{\alpha}$ using $\bar{\alpha}_{\bar{p}_1},\dots,\bar{\alpha}_{\bar{p}_k}$
		\State Compute $\bar{\epsilon}= \bar{\alpha}-\bar{f}+\bar{h}$		
	\end{algorithmic}
\end{algorithm}

The overall computational cost of the scheme is analyzed in the following.

\begin{proposition}\label{prop:complex}
	Assume that $\bar{p}_1,\dots,\bar{p}_k$ are all $v$-dimensional vectors and $v=m/k$. Then, the computational cost required by Algorithm \ref{alg:eccheo} is $	\mathcal{O}(mn^2)+\mathcal{O}\big(\frac{m^3+m^2nk+mn^2k^2}{k^2}\big)$.
\end{proposition} 
\begin{proof}
	By using, e.g., the singular value decomposition, a cost $\mathcal{O}(mn^2)$ is required for computing $\mathsf{G}$ and $\mathsf{L}$. Then, by taking into account matrix multiplications and (pseudo)inversions, for each $\ell=1,\dots,k$ we have a cost of $\mathcal{O}(v^3+v^2n+vn^2)$. Putting together, we get $		\mathcal{O}(mn^2)+k\mathcal{O}(v^3+v^2n+vn^2).$ Finally, by using $v=m/k$ the thesis follows.
\end{proof}
To compare what obtained to the computational cost of a standard CV implementation, which is outlined in Subsection \ref{sec:xy} and detailed in Algorithm \ref{alg:esatto}, note that when $m=n$ the result in Proposition \ref{prop:complex} becomes $\mathcal{O}(m^3)+\mathcal{O}\big(\frac{m^3+m^3k+m^3k^2}{k^2}\big)$. Therefore, if $k\approx m$, the Surrogate CV algorithm computational effort is $\mathcal{O}(m^3)$, more efficient than the $\mathcal{O}(m^4)$ required by the standard approach.

\begin{algorithm}[h!]
	\caption{Exact $k$-fold CV}\label{alg:esatto}
	\begin{algorithmic}[1]
		\Require $\mathsf{G}$, $\mathsf{L}$, $\bar{g}$, $\bar{h}$, $k$
		\Ensure $\bar{e}$
		\State Set $k$ vectors $\bar{p}_1,\dots,\bar{p}_k$ of distinct indices in $\{1,\dots,m\}$
		\For{$\ell=1,\dots,k$} 	
		\State $\bar{t}=\bigcup_{i\neq \ell}{\bar{p}_i}$ (training indices)
		\If{$m<n$}
		\State $\bar{w}=[m+1,\dots,n]$
		\State $\bar{\tau}=\bar{t}\cup \bar{w}$ (augmented training indices because of the added centers)
		\Else
		\State $\bar{\tau}=\bar{t}$
		\EndIf
		\State Compute $\bar{e}_{\bar{p}_\ell}= \bar{h}_{\bar{p}_\ell}-\mathsf{L}_{\bar{p}_\ell,\bar{\tau}}(\mathsf{G}_{\bar{t},\bar{\tau}})^{+}\bar{g}_{\bar{t}}$
		\EndFor
		\State Complete $\bar{e}$ using $\bar{e}_{\bar{p}_1},\dots,\bar{e}_{\bar{p}_k}$	
	\end{algorithmic}
\end{algorithm}

In our numerical tests, we will also compare our Surrogate CV to the \textit{empirical} Rippa-like strategy that has been employed in some previous works, which is detailed in Algorithm \ref{alg:ehnoeh}. Differently with respect to Algorithms \ref{alg:eccheo} and \ref{alg:esatto}, we point out that such an empirical scheme is employable in the LOOCV case only, and with $X=Y$. On the other hand, its related computational cost is limited to $\mathcal{O}(m^3)$.

\begin{algorithm}[h!]
	\caption{Empirical LOOCV}\label{alg:ehnoeh}
	\begin{algorithmic}[1]
		\Require $\mathsf{G}$, $\bar{g}$
		\Ensure $\bar{\eta}$
		\State Compute $\bar{c}=\mathsf{G}^{-1}\bar{g}$
		\State Compute $\bar{\eta}=\bar{c}./\mathrm{diag}(\mathsf{G}^{-1})$	
	\end{algorithmic}
\end{algorithm}

A Matlab implementation of the proposed Surrogate CV Algorithm \ref{alg:eccheo} is available at
\begin{equation*}
	\texttt{https://github.com/cesc14/RippaExtCV}\:.
\end{equation*}
\section{Numerics}\label{sec:numerics}

In the following, we test the discussed Surrogate CV algorithm by carrying out different numerical experiments that aim to show the benefits of the proposed scheme: restricted computational times and accurate approximation of the exact validation error. Therefore, we consider the following elliptic problem (see \cite[Section 39.1]{Fasshauer07}), which will be addressed via different validation strategies and in various settings.
\begin{equation*}
	\begin{split}
		& \Delta u(\bar{x}) =-\frac{5}{4}\pi^2\sin{(\pi x_1)}\cos{\bigg(\frac{\pi x_2}{2}\bigg)},\quad \bar{x}=(x_1,x_2)\in\Omega=[0,1]^2,\\
		&  u(\bar{x})= \sin(\pi x_1),\quad \bar{x}\in\Gamma_1,\\
		&  u(\bar{x})= 0,\quad \bar{x}\in\Gamma_2,
	\end{split}
\end{equation*}
where $\Gamma_1=\{\bar{x}\in\Omega\:|\: x_2=0\}$ and $\Gamma_2=\partial\Omega\setminus\Gamma_1$. The exact solution of this problem is
\begin{equation*}
	u(\bar{x})=\sin{(\pi x_1)}\cos{\bigg(\frac{\pi x_2}{2}\bigg)}.
\end{equation*}
Moreover, we will consider the following RBFs, that depend on a shape parameter that needs to be tuned,
\begin{equation*}
	\begin{split}
		& \varphi_{M,\varepsilon}(r)= e^{-\varepsilon r}(1+\varepsilon r) ,\quad\textit{Matérn $C^2$},\\
		& \varphi_{I,\varepsilon}(r)= \frac{1}{\sqrt{1+(\varepsilon r)^2}},\quad\textit{Inverse Multiquadrics $C^\infty$}.
	\end{split}
\end{equation*}
To construct the collocation schemes, we also consider
\begin{equation*}
	\begin{split}
		& \Delta\varphi_{M,\varepsilon}(r)= \varepsilon^2e^{-\varepsilon r}(\varepsilon r - 1),\\
		& \Delta\varphi_{I,\varepsilon}(r)= \frac{\varepsilon^2((\varepsilon r)^2-2)}{(1+(\varepsilon r)^2)^{5/2}},\quad \Delta\Delta\varphi_{I,\varepsilon}(r)= \frac{3\varepsilon^4(3(\varepsilon r)^4-24(\varepsilon r)^2+8)}{(1+(\varepsilon r)^2)^{9/2}}
	\end{split}
\end{equation*}

The experiments have been carried out in Matlab on a Intel(R) Core(TM) i7-1165G7 CPU@2.80GHz processor.

\subsection{Test 1: computational times and validation accuracy with Kansa's approach}\label{sec:kansa_comput}

Let $\bar{\mu}=(4^2,8^2,\dots,28^2,32^2)$. For each $\mu\in\bar{\mu}$, we perform the following test.
\begin{enumerate}
	\item
	We consider a set of internal collocation Halton points $H_\mu$ \cite{Halton60}, along with a set of $\sqrt{\mu}$ boundary collocation points $B_{\sqrt{\mu}}$ that are equispaced on the boundary of $\Omega$ (see Figure \ref{fig:1}). Therefore, we have $m=\mu + \sqrt{\mu}$ collocation points $X=H_\mu\cup B_{\sqrt{\mu}}$.
	
	\begin{figure}[h!]
		\centering
		\includegraphics[width=0.32\linewidth]{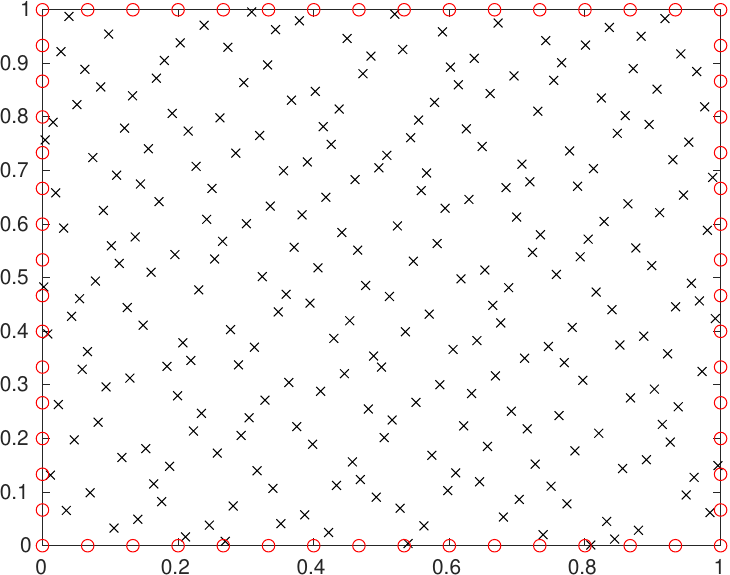}  
		\caption{The collocation Halton points (black crosses) and equispaced boundary points (red circles) in $\Omega$ for $\mu=256$.}
		\label{fig:1}
	\end{figure}
	\item 
	We let the centers coincide with the collocation points (case $X=Y$) and we implement Kansa's approach by setting	
	\begin{equation*}
		\lambda_{\bar{x}_j}(u) = \delta_{\bar{x}_j}(u)=u(\bar{x}_j)\quad j=1,\dots,m,
	\end{equation*}
	and
\begin{equation*}
	\begin{cases}
		\gamma_{\bar{x}_i}(u)=(\Delta u)_{|\bar{x}=\bar{x}_i}, & \quad \bar{x}_i\in\mathring{\Omega},\\
		\gamma_{\bar{x}_i}(u)=\delta_{\bar{x}_i}(u), & \quad \bar{x}_i\in\partial\Omega.	
	\end{cases}
\end{equation*}	
	
	\item 
	Let $\bar{\varepsilon}$ be a vector of $100$ shape parameter values between $2^{-5}$ and $2^{5}$ discretized in log-form. We want to choose $\varepsilon\in\bar{\varepsilon}$ so that the $L_2$-norm of the LOOCV error vector computed on $X$ is minimized.
	
	\item 
	Therefore, for each $\varepsilon\in\bar{\varepsilon}$ we compute the LOOCV ($k=m$) error via three different strategies.
	\begin{itemize}
		\item 
		Exact LOOCV: we compute the classical LOOCV error vector $\bar{e}$.
		\item 
		Surrogate LOOCV: we calculate an approximate value of the vector of LOOCV error $\bar{\epsilon}$ by using the proposed Algorithm \ref{alg:eccheo}.
		\item 
		Empirical LOOCV: we compute an \textit{inexact} LOOCV error vector $\bar{\eta}$ by using the empirical Algorithm \ref{alg:ehnoeh}.
		
	\end{itemize}

\end{enumerate}

\begin{figure}[h!]
	\centering
	\includegraphics[width=0.32\linewidth]{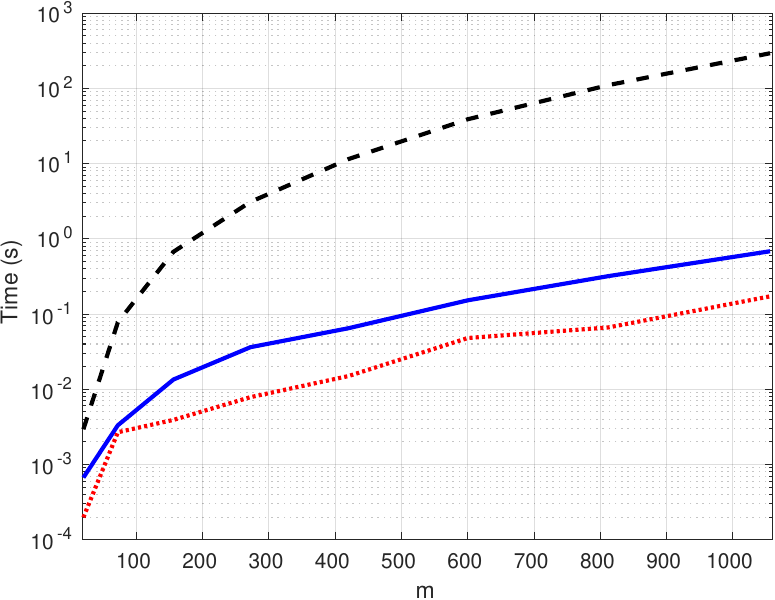}
	\includegraphics[width=0.32\linewidth]{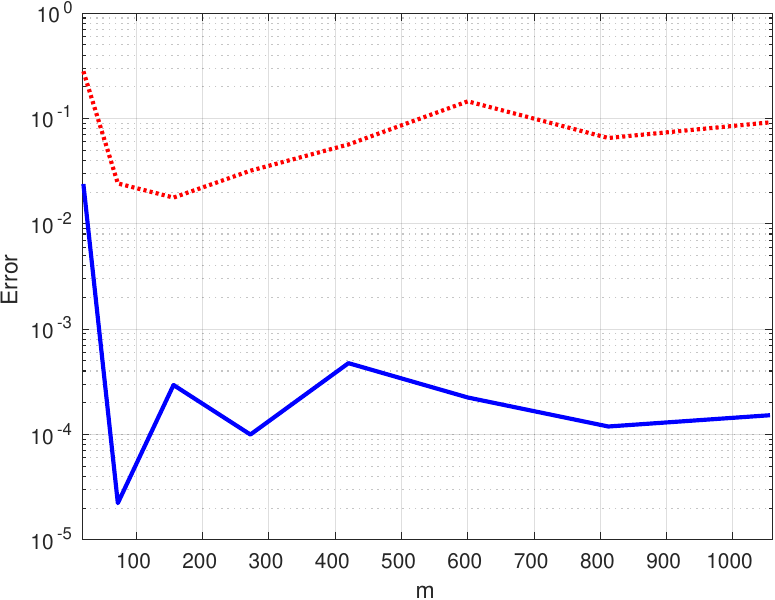}  
	\includegraphics[width=0.32\linewidth]{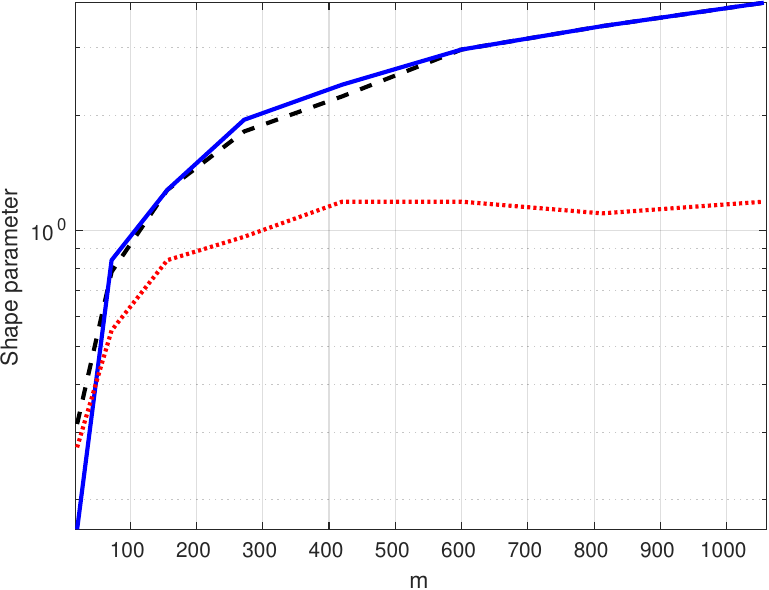}  
	\caption{RBF: $\varphi_{I,\varepsilon}$. On the left-hand side, the computational time employed varying $m$ by Exact LOOCV (dashed black), Surrogate LOOCV (solid blue) and Empirical LOOCV (dotted red) in validating the resulting best shape parameter value $\varepsilon^{\star}$, which is displayed on the right-hand side. On the center, we display $|\lVert \bar{e}\lVert_2-\lVert\bar{\epsilon}\lVert_2|$ (solid blue) and $|\lVert\bar{e}\lVert_2-\lVert\bar{\eta}\lVert_2|$ (dotted red) corresponding to $\varepsilon^{\star}$.}
	\label{fig:2}
\end{figure}

\begin{figure}[h!]
	\centering
	\includegraphics[width=0.32\linewidth]{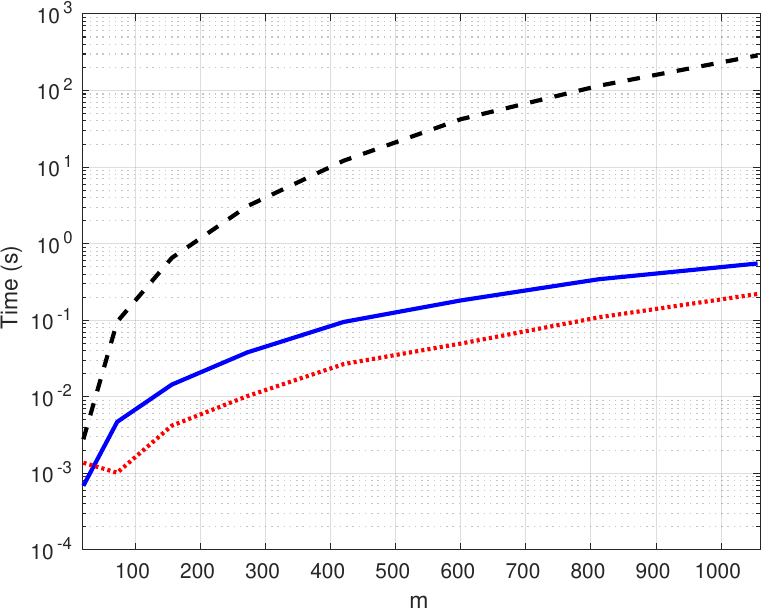}  
	\includegraphics[width=0.32\linewidth]{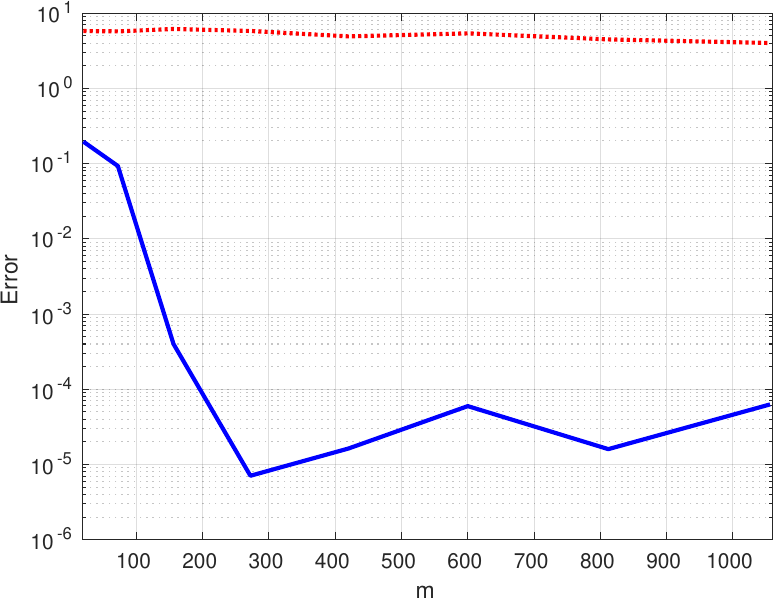}  
	\includegraphics[width=0.32\linewidth]{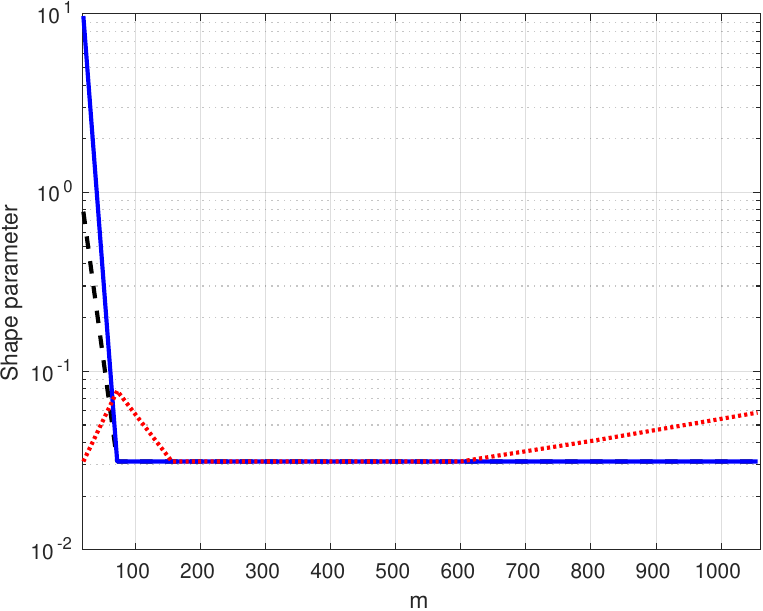}  
	\caption{RBF: $\varphi_{M,\varepsilon}$. On the left-hand side, the computational time employed varying $m$ by Exact LOOCV (dashed black), Surrogate LOOCV (solid blue) and Empirical LOOCV (dotted red) in validating the resulting best shape parameter value $\varepsilon^{\star}$, which is displayed on the right-hand side. On the center, we display $|\lVert \bar{e}\lVert_2-\lVert\bar{\epsilon}\lVert_2|$ (solid blue) and $|\lVert\bar{e}\lVert_2-\lVert\bar{\eta}\lVert_2|$ (dotted red) corresponding to $\varepsilon^{\star}$.}
	\label{fig:3}
\end{figure}

In Figure \ref{fig:2} we report the results obtained by considering $\varphi_{I,\varepsilon}$, while in Figure \ref{fig:3} we use as basis function $\varphi_{M,\varepsilon}$.

First, we observe that the proposed Surrogate LOOCV is definitely more efficient than Exact LOOCV, as expected, and slightly slower than Empirical LOOCV. Moreover, the central plots in Figures \ref{fig:2} and \ref{fig:3} show that Surrogate LOOCV computes a much more accurate validation error with respect to the empirical scheme. As a consequence, it is more likely for Surrogate LOOCV to select a shape parameter that coincides with the choice of the exact scheme, differently from Empirical LOOCV.

\subsection{Test 2: a numerical experiment with the Hermite method}

In the following, we set $\mu=256$ and the set of collocation points and centers $X=H_\mu\cup B_{\sqrt{\mu}}$. However, differently with respect to the previous section, here we restrict to $\varphi_{I,\varepsilon}$ and we test our algorithm in dealing with the Hermite method, which is characterized by the choice

\begin{equation*}
	\begin{cases}
		\lambda_{\bar{x}_i}(u)=\gamma_{\bar{x}_i}(u)=(\Delta u)_{|\bar{x}=\bar{x}_i}, & \quad \bar{x}_i\in\mathring{\Omega},\\
		\lambda_{\bar{x}_i}(u)=\gamma_{\bar{x}_i}(u)=\delta_{\bar{x}_i}(u), & \quad \bar{x}_i\in\partial\Omega.	
	\end{cases}
\end{equation*}	

Then, we perform steps 3. and 4. of Subsection \ref{sec:kansa_comput} and we focus on the obtained validation errors and best shape parameter values. The results are reported in Table \ref{tab:1}, where we highlight the $L_2$-norm of the obtained validation error vector related to the chosen shape parameter value.

\begin{table}[h!]

	\centering
		\begin{tabular}{ | c | c c |}
			\hline
			LOOCV & Best validation error & Best $\varepsilon$ \\
			\hline 
			Exact & $5.8778\textrm{E}-05$ & $1.5763$ \\
			\hline
			Surrogate & $9.2323\textrm{E}-06$ & $1.5763$ \\
			\hline		
			Empirical & $5.2208\textrm{E}-04$ & $0.6344$ \\
			\hline
		\end{tabular}
	\caption{Results obtained via the Hermite approach.
}\label{tab:1}
\end{table}

Also in this Hermite case, Surrogate LOOCV outperforms Empirical LOOCV in terms of accuracy.

\subsection{Test 3: dealing with a non-square collocation matrix}

\begin{figure}[h!]
	\centering
	\includegraphics[width=0.32\linewidth]{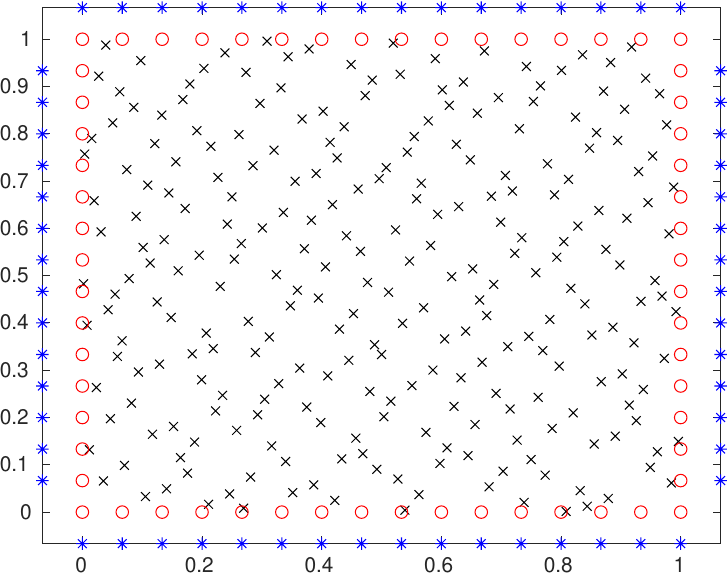}  
	\caption{The collocation Halton points (black crosses) and equispaced boundary points (red circles) in $\Omega$, the added centers outside $\Omega$ (blue stars) for $\mu=256$.}
	\label{fig:4}
\end{figure}

Here, we set again $\mu=256$ and $X=H_\mu\cup B_{\sqrt{\mu}}$, but we choose a different set of centers $Y$. Precisely, we add $\sqrt{\mu}$ centers that lie outside $\Omega$ as depicted in Figure \ref{fig:4}. This is a well-known strategy adopted to improve the accuracy of the collocation scheme nearby the boundary of the domain \cite[Section 39.1]{Fasshauer07}. Then, we consider Kansa's approach, with basis function $\varphi_{M,\varepsilon}$, and we proceed as in the previous subsection. However, note that in this case the matrices $\mathsf{G}$ and $\mathsf{L}$ are not squared, therefore the empirical LOOCV of Algorithm \ref{alg:ehnoeh} is not applicable.

\begin{table}[h!]
	\centering
	\begin{tabular}{ | c | c c |}
		\hline
			LOOCV & Best validation error & Best $\varepsilon$ \\
		\hline 
		Exact & $2.1539\textrm{E}+00$ & $0.0313$ \\
		\hline
		Surrogate & $2.1543\textrm{E}+00$ & $0.0313$ \\
		\hline		
	\end{tabular}
	\caption{Results obtained with a non-square collocation matrix.
	}\label{tab:2}
\end{table}

We observe that our proposed Surrogate LOOCV is fairly accurate in this non-square collocation setting too, being $\lVert \bar{\epsilon}\lVert_2$ close to $\lVert \bar{e}\lVert_2$.

\subsection{Test 4: varying $k$ in $k$-fold CV}

\begin{figure}[h!]
	\centering
	\includegraphics[width=0.32\linewidth]{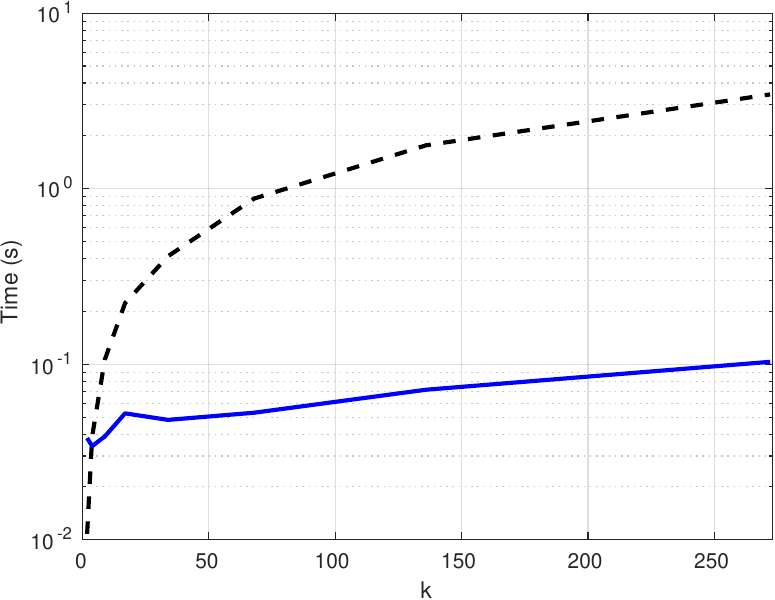}  
	\includegraphics[width=0.32\linewidth]{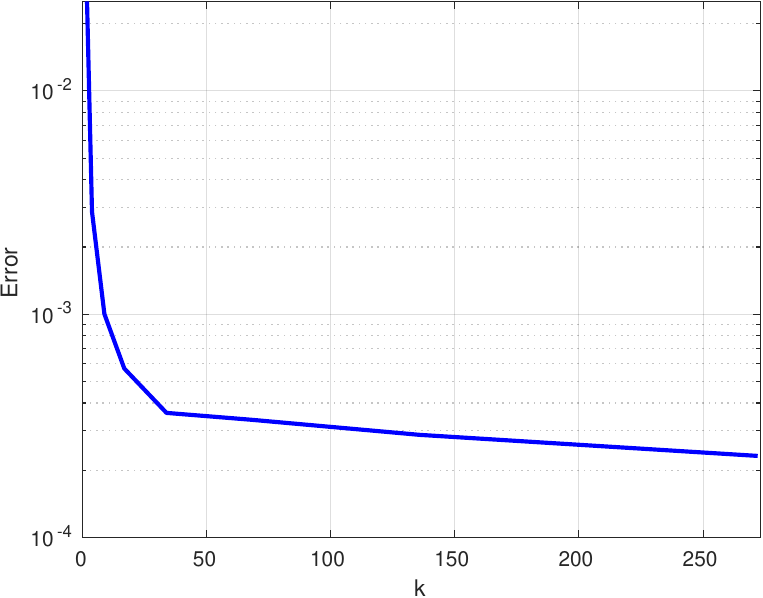}  
	\caption{On the left-hand side, the computational time employed varying $k$ by Exact CV (dashed black) and Surrogate CV (solid blue) in validating the resulting best shape parameter value. On the left, we display $|\lVert \bar{e}\lVert_2-\lVert \bar{\epsilon}\lVert_2|$.}
	\label{fig:5}
\end{figure}

In this subsection, our purpose is to analyze the behavior of the proposed Surrogate CV varying the number of folds $k$, in terms of both computational time and accuracy. To do this, we set again $\mu=256$ and $X=Y=H_\mu\cup B_{\sqrt{\mu}}$ and employ Kansa's approach with basis function $\varphi_{M,\varepsilon}$. However, here steps 3. and 4. of Subsection \ref{sec:kansa_comput} are repeated for $k\in\{\lfloor m/2^i\rfloor\:|\:i=0,\dots,7\}$. We report the achieved results in Figure \ref{fig:5}.

As far as the computational times are concerned, as expected Exact CV is faster than Surrogate CV for very small values of $k$ only. About the accuracy of the proposed scheme varying $k$, the right-hand plot in Figure \ref{fig:5} provides an interesting insight: the approximation of the exact validation error gets worser as $k$ becomes smaller. This fact can be related to \eqref{eq:lss}, from which the inexactness of our Surrogate CV originates. Indeed, a small value of $k$ implies a large number $v$ of validation samples, and therefore the model is built on a restricted number of data. Consequently, it is natural to deviate from the RBF-PS solution. Nevertheless, we observe that the scheme can retain a \textit{good} accuracy for any $k$.

\section{Conclusions}\label{sec:conclusions}

In this work, we proposed a new surrogate $k$-fold CV scheme for the RBF collocation setting, which is inspired by the extended Rippa's algorithm employed in the interpolation framework. The proposed method is more efficient than a straightforward calculation of exact $k$-fold CV, especially when $k$ is \textit{large}, and can also be used when additional centers besides collocation points are taken into consideration. Moreover, in the case of LOOCV, it provides a more accurate approximation of the exact validation error with respect to the Rippa-like empirical LOOCV approach that is used in many previous works. The discussed Surrogate CV Algorithm \ref{alg:eccheo} has been first analyzed in details, and then tested in various numerical experiments. Future work consists in studying the proposed scheme in other collocation settings, also with conditionally positive definite kernels, and in evaluating possible modifications to further enhance its efficiency \cite{Ling22}.


\begin{thebibliography}{99}
	
	\bibitem{Campagna20} 
	\textsc{R. Campagna, S. Cuomo, S. De Marchi, E. Perracchione, G. Severino}, 
	\emph{A stable meshfree PDE solver for source-type flows in porous media}, Appl. Numer. Math \textbf{149} (2020), pp. 30--42.
	
	\bibitem{Cavoretto22} 
	\textsc{R. Cavoretto}, 
	\emph{Adaptive LOOCV-based kernel methods for solving time-dependent BVPs}, Appl. Math. Lett. \textbf{429} (2022), 127228.
	
	\bibitem{Cavoretto20a}
	\textsc{R. Cavoretto, A. De Rossi},
	\emph{A two-stage adaptive scheme based on RBF collocation for solving elliptic PDEs}, Comput. Math. Appl. \textbf{79} (2020), pp. 3206-3222.
	
	\bibitem{Cavoretto20b}
	\textsc{R. Cavoretto, A. De Rossi},
	\emph{An adaptive LOOCV-based refinement scheme for RBF collocation methods over irregular domains}, Appl. Math. Lett. \textbf{103} (2020), 106178.
	
	\bibitem{Cavoretto18}
	\textsc{R. Cavoretto, A. De Rossi, E. Perracchione}, \emph{Optimal Selection of Local Approximants in RBF-PU Interpolation}, J Sci Comput \textbf{74} (2018), pp. 1--22.
	
	\bibitem{Chen20}
	\textsc{C. S. Chen, A. Karageorghis, H. Zheng}, \emph{Improved RBF Collocation Methods for Fourth Order Boundary Value Problems}, Commun. Comput. Phys. \textbf{27} (2020), pp. 1530--1549.
	
	\bibitem{Chen22}
	\textsc{Y.-T. Chen, C. Li, L.-Q. Yao, Y. Cao}, \emph{A Hybrid RBF Collocation Method and Its Application in the Elastostatic Symmetric Problems}, Symmetry \textbf{14}(7) (2022), 1476.

	\bibitem{Chen22a}
	\textsc{M. Chen, L. Ling, Y. Su}, \emph{Solving interpolation problems on surfaces stochastically and greedily}, Dolomites Res. Notes Approx. \textbf{15}(3) (2022), pp. 26--36.	

	\bibitem{Chiappa20}
	\textsc{A. Chiappa, C. Groth, A. Reali, M. Evangelos Biancolini}, \emph{A stress recovery procedure for laminated composite plates based on strong-form equilibrium enforced via the RBF Kansa method}, Compos. Struct. \textbf{244} (2020), 112292.
	
	\bibitem{Chiu20}
	\textsc{S.N. Chiu, L. Ling, M. McCourt}, \emph{On variable and random shape Gaussian interpolations}, Appl. Math. Comput. \textbf{377} (2020), 125159.
	
	\bibitem{Chu14}
	\textsc{F. Chu, L. Wang, Z. Zhong, J. He}, \emph{Hermite radial basis collocation method for vibration of functionally graded plates with in-plane material inhomogeneity}, Comput. Struct. \textbf{142} (2014), pp. 79--89.		
	
	\bibitem{Dehghan15}
	\textsc{M. Dehghan, M. Abbaszadeh, A. Mohebbi}, \emph{An implicit RBF meshless approach for solving the time fractional nonlinear sine-Gordon and Klein–Gordon equations}, Eng. Anal. Bound. Elem. \textbf{50} (2015), pp. 412--434.
	
	\bibitem{Dehghan14}
	\textsc{M. Dehghan, V. Mohammadi}, \emph{The numerical solution of Fokker-Planck equation with radial basis functions (RBFs) based on the meshless technique of Kansa's approach and Galerkin method}, Eng. Anal. Bound. Elem. \textbf{47} (2014), pp. 38--63.
	
	\bibitem{Fasshauer07}
	\textsc{G.E. Fasshauer},
	\emph{Meshfree Approximations Methods with \textsc{Matlab}}, 
	World Scientific, Singapore, 2007.
	
	\bibitem{Fasshauer15}
	\textsc{G.E. Fasshauer, M.J. McCourt}, 
	\emph{Kernel-based Approximation Methods Using} \textsc{Matlab}, 
	World Scientific, Singapore, 2015.
	
	\bibitem{Fasshauer07a}
	\textsc{G.E. Fasshauer, J.G. Zhang}, \emph{On choosing \lq\lq optimal\rq\rq shape parameters for RBF approximation}, Numer. Algorithms \textbf{45} (2007), pp. 345--368.
	
	\bibitem{Fornberg07} 
	\textsc{B. Fornberg, J. Zuev}, \emph{The Runge phenomenon and spatially variable shape parameters in RBF interpolation}, Comput. Math. Appl. \textbf{54}(3) (2007), pp. 379--398.
	
	\bibitem{Gherlone12} 
	\textsc{M. Gherlone, L. Iurlaro, M. Di Sciuva}, \emph{A novel algorithm for shape parameter selection in radial basis functions collocation method}, Compos. Struct. \textbf{94}(2) (2012), pp. 453--461.
	
	\bibitem{Golbabai15}
	\textsc{A. Golbabai, E. Mohebianfar, H. Rabiei}, \emph{On the new variable shape parameter strategies for radial basis functions}, J. Comput. Appl. Math. \textbf{34} (2015), pp. 691--704.
	
	\bibitem{Golub79}
	\textsc{G. H. Golub, M. Heath, G. Wahba}, \emph{Generalized cross-validation as a method
	for choosing a good ridge parameter}, Technometrics \textbf{21}(2) (1979), pp. 215--223.

	\bibitem{Halton60}
	\textsc{J. H. Halton}, \emph{On the efficiency of certain quasi-random sequences of points in evaluating multi-dimensional integrals}, Numer. Math. \textbf{2} (1960), pp. 84--90.

	\bibitem{Hon01}
	\textsc{Y. C. Hon, R. Schaback}, \emph{On nonsymmetric collocation by radial basis functions}, Appl. Math. Comput. \textbf{119} (2001), pp. 177--186.

	\bibitem{Kansa90} 
	\textsc{E. J. Kansa}, 
	\emph{Multiquadrics–A scattered data approximation scheme with applications to computational fluid-dynamics–II solutions to parabolic, hyperbolic and elliptic partial differential equations}, Comput. Math. Appl. \textbf{19} (1990), pp. 147--161.

	\bibitem{Kansa92} 
	\textsc{E. J. Kansa, R. Carlson}, 
	\emph{Improved accuracy of multi-quadric interpolation using variable shape parameters}, Comput. Math. Appl. \textbf{24} (1992), pp. 20--99.
	
	\bibitem{Karageorghis21} 
	\textsc{A. Karageorghis, D. Tappoura, C.S. Chen}, 
	\emph{The Kansa RBF method with auxiliary boundary centres for fourth order boundary value problems}, Math. Comput. Simul. \textbf{181} (2021), pp. 581--597.
	
	\bibitem{Katsiamis20} 
	\textsc{A. Katsiamis, A. Karageorghis}, 
	\emph{Kansa radial basis function method with fictitious centres for solving nonlinear boundary value problems}, Eng. Anal. Bound. Elem. \textbf{119} (2020), pp. 293--301.

	\bibitem{Kazem12} 
	\textsc{S. Kazem, J.A. Rad, K. Parand}, 
	\emph{Radial basis functions methods for solving Fokker–Planck equation}, Eng. Anal. Bound. Elem. \textbf{36}(2) (2012), pp. 181--189.
	
	\bibitem{Krowiak19} 
	\textsc{A. Krowiak, J. Podgórski}, 
	\emph{On choosing a value of shape parameter in radial basis function collocation methods}, AIP Conference Proceedings \textbf{2116}(1) (2019), 450020.
	
	\bibitem{LaRocca05} 
	\textsc{A. La Rocca, A. Hernandez Rosales, H. Power}, 
	\emph{Radial basis function Hermite collocation approach for the solution of time dependent convection–diffusion problems}, Eng. Anal. Bound. Elem. \textbf{29}(4) (2005), pp. 359--370.
		
	\bibitem{LaRocca06} 
	\textsc{A. La Rocca, H. Power}, 
	\emph{A Hermite radial basis function collocation approach for the numerical simulation of crystallization processes in a channel}, Commun. Numer. Meth. Engng. \textbf{22} (2006), pp. 119--135.
	
	\bibitem{Ling22} 
	\textsc{L. Ling, F. Marchetti}, 
	\emph{A stochastic extended Rippa’s algorithm for LpOCV}, Appl. Math. Lett. \textbf{129} (2022), 107955.
	
	\bibitem{Liu15} 
	\textsc{XY. Liu, A. Karageorghis, C.S. Chen}, 
	\emph{A Kansa-Radial Basis Function Method for Elliptic Boundary Value Problems in Annular Domains}, J. Sci. Comput. \textbf{65} (2015), pp. 1240--1269.

	\bibitem{Ma21} 
	\textsc{X. Ma, B. Zhou, S. Xue}, 
	\emph{A meshless Hermite weighted least-square method for piezoelectric structures}, Appl. Math. Comput. \textbf{400} (2021), 126073.

	\bibitem{Marchetti21} 
	\textsc{F. Marchetti}, 
	\emph{The extension of Rippa’s algorithm beyond LOOCV}, Appl. Math. Lett. \textbf{120} (2021), 107262.
	
	\bibitem{Mongillo11}
	\textsc{M. Mongillo}, \emph{Choosing Basis Functions and Shape Parameters for Radial Basis Function Methods}, SIAM SIURO publications \textbf{4} (2011).
	
	\bibitem{Rippa99}
	\textsc{S. Rippa}, \emph{An algorithm for selecting a good value for the parameter $c$ in radial basis function interpolation}, Adv. Comput. Math. \textbf{11} (1999), pp. 193--210.
	
	\bibitem{Roque10}
	\textsc{C.M.C. Roque, A.J.M. Ferreira}, \emph{Numerical Experiments on Optimal ShapeParameters for Radial Basis Functions}, Numer. Meth. Partial Diff. Eqs. \textbf{26} (2010), pp. 675--689.
	
	\bibitem{Schaback07}
	\textsc{R. Schaback}, \emph{Convergence of Unsymmetric Kernel-Based Meshless Collocation Methods}, SIAM J Numer Anal \textbf{45}(1) (2007), pp. 333--351.	
	
	\bibitem{Scheuerer11}
	\textsc{M. Scheuerer}, \emph{An alternative procedure for selecting a good value for the parameter c in RBF-interpolation}, Adv. Comput. Math. \textbf{34} (2011), pp. 105--126.
	
	\bibitem{Trahan03}
	\textsc{C.J. Trahan, R.E. Wyatt}, \emph{Radial basis function interpolation in the quantum trajectory method: optimization of the multi-quadric shape parameter}, J. Comput. Phys. \textbf{185} (2003), pp. 27--49.	
	
	\bibitem{Uddin14}
	\textsc{M. Uddin}, \emph{On the selection of a good value of shape parameter in solving time-dependent partial differential equations using RBF approximation method}, Appl. Math. Model. \textbf{38} (2014), pp. 135--144.	
	
	\bibitem{Yang18}
	\textsc{F. Yang, L. Yan, L. Ling}, \emph{Doubly stochastic radial basis function methods}, J. Comput. Phys. \textbf{363} (2018), pp. 87--97.
	
	\bibitem{Wang18}
	\textsc{F. Wang, W. Chen, C. Zhang, Q. Hua}, \emph{Kansa method based on the Hausdorff fractal distance for Hausdorff derivative Poisson equations}, Fractals \textbf{26}(4) (2018), 1850084.
	
	\bibitem{Wendland05}
	\textsc{H. Wendland}, 
	\emph{Scattered Data Approximation}, 
	Cambridge Monogr. Appl. Comput. Math., vol. 17, Cambridge Univ. Press, Cambridge, 2005.
	
\end{thebibliography}

\end{document}